\documentstyle[12pt,amscd, amssymb]{amsart}

\newtheorem{theo}{THEOREM}[section]
\newtheorem{lemma}[theo]{Lemma}
\newtheorem{cor}[theo]{Corollary}
\newtheorem{prop}[theo]{Proposition}
\theoremstyle{remark}

\newtheorem{rem}[theo]{Remark}
\newtheorem{case}{Case} 
\newcommand{\brref}[1]{(\ref{#1})}
\newcommand{\calo}{\cal {O}}
\newcommand{\tensor}{\otimes}
\newcommand{\map}[3]{$#1$~:~$#2~\rightarrow~#3$}
\newcommand{\Proj}[1]{{\Bbb P}(#1)}

\newcommand{\restrict}[2]{{#1}_{\mid _{#2}}}
\newcommand{\oof}[2]{{\cal O}_{#1}({#2})}
\newcommand{\oofp}[2]{{\cal O}_{{\Bbb P}^{#1}}({#2})}
\newcommand{\iof}[2]{{\cal I}_{#1}({#2})}
\newcommand{\xel}{(X, L)}
\newcommand{\sel}{(S, L)}
\newcommand{\elonebar}{\overline{\ell}_1}
\newcommand{\eltwobar}{\overline{\ell}_2}
\newcommand{\Pin}[1]{{\Bbb P}^{#1}}
\newcommand{\taut}[1]{{\cal O}_{{\Bbb P}(#1)}(1)}

\newcommand{\lra}{\longrightarrow}
\newcommand{\G}{\Gamma}
\newcommand{\Gs}{\Gamma^*}
\title[Projective Normality]{On the Projective Normality of Smooth
Surfaces of
Degree Nine}
\author[G. Besana ]{Gian Mario Besana$^*$}
\author[S. Di Rocco] {Sandra Di Rocco$^{**}$}
\address[G. Besana]{\small Department of Mathematics - Eastern Michigan University - Ypsilanti MI- U.S.A.48197} 
\address[S. Di Rocco]{\small Department of Mathematics- KTH - Royal Institute of Technology -
S 100 44 Stockholm - SWEDEN} 
\email[G. Besana] {gbesana@@emunix.emich.edu}
\email[S. Di Rocco]{sandra@@math.kth.se } 
\thanks{$^*$ Partially supported by M.U.R.S.T.  $^{**}$ Partially supported
by C.N.R.)}

\begin{document}
\maketitle

\begin{abstract}
The projective normality of smooth, linearly normal   surfaces
of degree 9 in $\Pin{N}$ is studied. All non projectively normal
surfaces which 
are not scrolls over a curve are classified. Results on the projective
normality of
surface scrolls are also given.
\end{abstract}

\section{introduction}
Smooth projective varieties of small degree have been classified over
the years
and thoroughly studied, e.g. \cite{Io1}, \cite{Io2}, \cite{ok2},
\cite{ok8},
\cite {Ale}, \cite{ADS}.  A variety $X~\subset~\Pin{n}$ is {\it
projectively normal} if
the 
maps
$H^0(\Pin{n},\oofp{n}{k}) \to H^0(X, \oof{X}{k})$ are surjective for all
$k
\ge 1$ or
in other words if hypersurfaces of degree $k$ cut complete linear
systems
on $X$ for every $k\ge 1.$ In
\cite{Alibaba} the projective normality of varieties of degree $d \le 8$
of
any dimension  was
investigated. This work is concerned with the projective normality of
smooth
projective surfaces, embedded by the complete linear system associated
with a
very ample line bundle
$L$   of degree
$d=9.$ Such surfaces are either embedded in $\Pin{4}$ or have sectional
genus $g \le 7.$ Therefore they are completely classified in
\cite{au-ra} and
\cite{LiAq}.
One of the reasons that brought us to look at this question is our
desire
to find examples for a long standing problem in adjunction theory.
 Andreatta
\cite{ce}, followed by a generalization by Ein and Lazarsfeld
\cite{EL},
posed the problem of classifying smooth $n$-dimensional varieties $(X,
\cal{L})$
polarized
with a very ample line
bundle $\cal{L}$, such that the adjoint linear system $|L| = |K + (n-1)
\cal{L}|$ gives an
embedding
which is not projectively normal. Andreatta and Ballico \cite{an-ba1}
gave examples of
surfaces $(S, \cal{L})$ with the above behavior, where $d=\deg{S} = 10$
under the
adjoint embedding.  Alzati, Bertolini and the first  author  in
\cite{Alibaba} found no
example with
$d\le 8.$  After a  detailed check of the non projectively normal
surfaces
found in this work no examples were found except possibly a blow
up of an elliptic $\Pin{1}$-bundle whose existence is uncertain. See
section  
\ref{K+L} for details.

 Our findings concerning the  projective normality of surfaces of degree
nine are
collected in the  the following theorem   (see
\brref{notation} for notation):

\begin{theo}
\label{thetheorem} Let $S$ be a smooth  surface embedded by the complete
linear system associated with a very ample line
bundle
$L$  as a surface  of degree $9$ and sectional genus $g$ in $\Pin{N}. $
Assume $\sel$
is not a scroll over a curve. Then
$(S, L)$ fails to be projectively normal if and only if it belongs to
the
following list:
\begin{center}
\begin{tabular}{|c|c|c|c|} \hline
$\Pin{N} $ & $g$& $S$& $L$\\ \hline\hline
$\Pin{5}$ & $4$ & $Bl_3X$ where $X$ is a $\Pin{1}$-bundle  & \\
 & & over an elliptic
curve,
$e=0$&
$2\frak{C_0}+3\frak{f}-\sum_iE_i$
\\
\hline
$\Pin{5}$ & $5$ & Rational conic bundle $S=Bl_{15} \bold{F_e}$ , $0\le e
\le 5$  &
$2 \frak{C_0} + (6+e)\frak{f} -
\sum_iE_i $
\\
\hline
$\Pin{4}$ & $6$  &$Bl_{10}(\Pin{2})$&$3p^*(\oofp{2}{1})-\sum_i4E_i$ \\
\hline
$\Pin{4}$ & $7$  &$Bl_{15}(\Pin{2})$&$9p^*(\oofp{2}{1})-\sum_1^63E_i-$\\
& & &$\sum_7^92E_j-\sum_{10}^{15}E_k $\\ \hline
$\Pin{4}$ & $6$ & Projection of an Enriques surface & \\
 &    &   of degree 10 in ${\bf P}^5$& cf. \cite{au-ra} \\ \hline
$\Pin{4}$ & $7$ & Minimal elliptic surface & cf. \cite{au-ra} \\ \hline
$\Pin{4}$ & $8$ & Minimal surface of general type& cf. \cite{au-ra}\\
\hline
\end{tabular}
\end{center}
\end{theo}
The projective normality of surfaces which are scrolls over a curve of
genus $g,$
not included in the above theorem, was also investigated. Results are
collected in
Proposition  \ref{scrollprop}. We were not able to prove or disprove the
projective
normality of  scrolls over trigonal curves of genus $3,4,5.$

The authors would like to thank Ciro Ciliberto, Antonio Lanteri and
Andrew J.
Sommese  for many helpful conversations and  W. Chach\'{o}lski for his
insight  and
patience.

\section{Background material}
\subsection{NOTATION}
\label{notation}
Throughout this article $S$ denotes  a smooth connected projective
surface
 defined over the complex field {\bf C}. Its
structure sheaf is denoted  by
${\cal O_S}$ and the canonical sheaf of holomorphic $2$-forms on $S$ is
denoted by
$K_S$.

For any coherent sheaf $\Im$ on $S$,  $h^i(\Im )$ is the complex
dimension of
$H^i(S,\Im)$ and
$\chi=\chi(\cal{O}_S)=\chi(S)=\sum_i(-1)^ih^i(\cal{O}_S).$
Let $L$ be a  line bundle on $S.$ If $L$ is ample the pair $\sel$ is
called
a {\it
polarized surface}.

The following notation is used:\\
$|L|$, the complete linear system associated with L;\\
$d = L^2,$ the degree of $L$;\\
$g=g(S, L)$, the {\it sectional genus} of $\sel$, defined by
$2g-2=L\cdot
(K_S+L).$ If $C\in |L|$ is an irreducible and reduced element then
$g=g(C)$
is the
arithmetic genus of $C$;\\
$\Delta \sel = \Delta = 2+L^2-h^0(L)$ the Delta genus of $\sel$;\\
 $\bold{F_e}$, the Hirzebruch surface of invariant $e$ ;\\
$E^*$ the dual of a vector bundle $E$.\\
Cartier divisors, their associated line bundles and the invertible
sheaves
of their
holomorphic sections are used with no distinction. Given two divisors
$L$
and $M$
we denote linear equivalence by
$L \sim M$ and numerical equivalence by $L \equiv M.$
The blow up of a surface $X$ at $n$ points is denoted by $p: S=Bl_nX \to
X.$
When $X$ is a
$\Pin{1}$-bundle over a curve with fundamental section $C_0$ and generic
 fibre $f$ it is
$Num(X) ={\Bbb Z}[C_0]
\oplus
{\Bbb Z}[f]$ and
 the following shorthand is used: $\frak{C_0}=p^*(C_0)$ and $\frak{f}=
p^*(f).$

 A polarized surface $(S, L)$
is a {\em scroll}
or a {\em conic bundle} over a curve $C$ if there exists a surjective
morphism $p: S
\to C$ with connected fibers and an ample line bundle
$H$ on
$C$ such that, respectively, $ K_S + 2L = p^*(H)$ or $ K_S + L =
p^*(H).$
If $\sel$ is a
scroll then $S$ is a $\Pin{1}$-bundle over $C$ and $L \cdot f = 1$ for
every fibre $f.$

 In section
\ref{genere6} the notion of {\it reduction} of a smooth polarized
surface
is shortly
used. The best  reference is
\cite{BESO}.
\subsection{\small CASTELNUOVO BOUND}
Let $C\subset \Pin{N}$, then by Castelnuovo's lemma
\begin{equation}
\label{cast}
g(C)\leq\left[ \frac{d-2}{N-2}\right ](d-N+1-
(\left[ \frac{d-2}{N-2}\right ]-1)\frac{N-2}{2})
\end{equation}
where$[x]$ denotes the greatest integer $\leq x$.

\subsection{\small PROJECTIVE NORMALITY}
Let $S$ be a surface embedded in $\Pin{N}.$ $S$ is said to be {\it
k-normal} if the map
$$H^0({\cal O}_{\Pin{N}}(k))\longrightarrow H^0({\cal O}_S(k))$$ is
surjective .
 $S$ is said to be {\it projectively normal} if it is $ k$-normal for
every $k\geq 1$.
An ample line bundle $L$ on $S$ is {\it normally generated} if
$S^kH^0(L)\to H^0(L^k)$ is surjective for every $k\geq 1$. If $L$ is
normally generated then it is very ample and $S$, embedded in $\Pin{N}$
via
$|L|$ is projectively normal. A polarized surface $(S, L)$ is said
to be projectively normal if $L$ is very ample and $S$ is
projectively normal under the embedding given by $L.$
A polarized surface $\sel$ has a {\it ladder} if there exists  an
irreducible and
reduced element
$C \in |L|.$ The ladder is said to be
{\it regular} if $H^0(S, L)\to H^0(C, \restrict{L}{C})$ is onto.   If
$L$ is  very ample  $\sel$ clearly  has a ladder.

We recall the following general result due to Fujita:

\begin{theo}[\cite{fu}]
\label{fujitatheo}
Let $(S, L)$ be a polarized surface with a ladder. Assume
$g(L)\geq\Delta$. Then \\ i) The ladder is regular if $d\geq 2\Delta
-1$;\\
ii) L is normally generated, $g=\Delta$ and $H^1(S, tL)= 0$ for any t,
if $d\geq 2\Delta +1$.
\end{theo}
    \subsection{\small $k$ - REGULARITY}
\label{kreg}
 A good reference is \cite{mu1}.
A coherent sheaf ${\cal F}$ over $\Pin{n}$ is  {\em k-regular} if
$h^i({\cal F}(k-i))=0$ for
all $i >0.$ If ${\cal F}$  is $k$-regular then it is $k+1$-regular. If
$X\subset  \Pin{n} $ is an irreducible variety such that ${\cal I}(X)$ is
$k$-regular
then the homogeneous ideal $I_X=\oplus H^0(\iof{X}{t})$ is generated in
degree $\le
k.$  This fact implies that if ${\cal I}_X$ is $k$-regular then  $X$
cannot
be embedded
with a $t\ge( k+1)$- secant line.
    \subsection{\small CLIFFORD INDEX}
 Good references are \cite{mart}, \cite{GL}.

 Let $C$ be a projective curve and
$H$ be any line bundle on $C$. The Clifford index of $H$  is defined as
follows:
$$cl(H)=d-2(h^0(H)-1).$$
The Clifford index of the curve is  $cl(C)=\text {min}\{cl(H) |
h^0(H)\geq 2 \text{ and }h^1(H)\geq 2 \}$.
$ H$ {\it contributes }to the Clifford index of $C$ if $h^0(H)\geq 2
\text{ and }h^1(H)\geq 2$ and $H$ {\it computes}
the  Clifford index of $C$ if $cl(C)=cl(L)$. For a general curve $C$
it is $cl(C)=\left [\frac{g-1}{2}\right ]$
 and in any case $cl(C)\leq\left [\frac{g-1}{2}\right ]$. By
Clifford's  theorem a  special line bundle $L$ on $C$ has $cl(L)\geq
0$ and
the equality holds if and only if $C$ is hyperelliptic and $L$ is a
multiple of the unique
$g^1_2$.\\ If
$cl(C)=1$ then $C$ is  either a plane quintic curve or a trigonal curve.
The following results dealing with the projective normality of curves
and relating it to the Clifford index are listed for the convenience of
the reader.
\begin{theo}[\cite{GL}]
\label{glcliff}
Let L be a very ample line bundle on a smooth irreducible complex
projective
curve C with:
$$deg(L)\geq 2g+1-2h^1(L)-cl(C)$$
then $L$ is normally generated.
\end{theo}

In the case of hyperelliptic curves, because  there are no special
very-ample line
bundles,  the following is true.
\begin{prop}[\cite{la-ma}] \label{hyper}
A hyperelliptic curve of genus $g$ has no normally generated line
bundles
of degree
$\leq 2g$.
\end{prop}
\begin{lemma}[ \cite{la-ma}]
\label{2norm}
Let $L$ be a base point free line bundle of degree $d\geq g+1$ on a
curve
of genus $g$. Then $L$ is normally generated if and only if the natural
map $H^0(L)\tensor H^0(L)\to H^0(2L)$ is onto.
\end{lemma}
It follows that  a polarized curve $(C, L)$ with
$d\geq g+1$ and $L$ very ample is  projectively normal if and only if it
is 2-normal.

The projective normality of a polarized surface $\sel$ will be often
established by investigating the   property
for a general hyperplane section. The main tools used are the following
results.
\begin{theo}[\cite{fu}]
\label{fujitatheo2}
Let $\sel\supset
(C,\restrict{L}{C})$ be a polarized surface with a ladder. If
the ladder is regular and $\restrict{L}{C}$ is normally generated then
$L$ is
 normally generated.
\end{theo}
\begin{lemma}[\cite{Alibaba}]
\label{besanaignorans}
Let $\sel\supset
(C,\restrict{L}{C})$ be a polarized surface with a regular ladder. 
Assume $h^1(L)=0$
and $\Delta=g.$
 Then $L$ is normally generated if and
only if $\restrict{L}{C}$ is normally generated.
\end{lemma}
\subsection{\small SURFACES EMBEDDED IN QUADRIC CONES}
\label{qcones}
As Lemma \ref{2norm} suggests, the hyperplane section technique
will often reduce the projective normality of a surface to its
$2$-normality. It is useful then to recall the detailed investigation
of surfaces in $\Pin{5}$ contained in singular quadrics, done in
\cite{Alibabaquad}.
  Let $\Gamma$ be a four dimensional quadric cone in
$\Pin{5}$ and let $\sigma:\Gs\lra\G$ be the blow up of $\G$ along the
vertex, with exceptional divisor $T$. Suppose $S\subset\G$ is a smooth
surface  and let
$S'$ be the strict transform of $S$ in $\Gs$ under $\sigma$.

The Chern classes of $S'$ and $\Gs$ satisfy the following
standard relation:
\begin{equation}
\label{DPF}
\restrict{c_2(\Gs)}{S'}=S'S'+\restrict{c_1(\Gs)}{{S'}}c_1(S')-
K_{S'}^2+c_2(S')
\end{equation}

If $rank(\G)=5$, \ $\G$ is a cone with vertex a point $P$ over a smooth
quadric $Q' \subset \Pin{4}.$
Following \cite{Alibabaquad}  let
 \begin{xalignat}{2}
C(W):&=&  \text{ the cone over the cycle $W \subset Q'$ with vertex $V$}
&
\notag
\\
\sigma:& \Gamma^* \to \Gamma  &\text{ the blow up map}& \notag \\
H_{Q'}: &=&   \text {the hyperplane section of $Q'$}&  \notag \\
l_{Q'}:&=&  \text {the generator of $A_1(Q')$}& \notag \\
p_{Q'}:&=& \text {the generator of $A_0(Q')$ }&\notag
\end{xalignat}
According to the above notation it is
\begin{gather}
Pic(Q') = <H_{Q'}> \notag \\
H^2_{Q'} = 2 l_{Q'} \notag\\
A_0(Q') =<p_{Q'}> \notag
\end{gather}
Further, let
\begin{gather}
Z =  \sigma^{-1}(C(H_{Q'})) \notag \\
S = \sigma^{-1}(C( l_{Q'})) \notag \\
F = \sigma^{-1}(C(p_{Q'})), \notag
\end{gather}
and denote by $\overline{H}$ the cycle $H_{Q'}$ in $\Gamma^*$ and
by $\overline{l}$ the
cycle $l_{Q'}$ in $\Gamma^* .$
The Chow Rings of $\Gamma^*$ are then given by:
\begin{gather}
Pic(\Gamma^*) = <Z, \tau>\\
A_2(\Gamma^*) = <\overline{H}, S>\\
A_1(\Gamma^*) = < \overline{l}, F>.
\end{gather}
and it is
\begin{xalignat}{3}
c_1(\Gamma^*) &= 2(\tau + Z) & c_2(\Gamma^*) & = Z^2 + 6\tau Z
&T&=\tau-Z
\end{xalignat}
\begin{xalignat}{2}
 T&=\tau-Z& S'&=\alpha \overline{H} +\beta X
\end{xalignat}
The intersection table is then the following:
\vskip .5 cm
\setlength{\unitlength}{ .4cm}
\begin{center}
\begin{picture}(13,13)
\put(1,0){\line(0,1){13}}
\put(5,0){\line(0,1){13}}
\put(9,0){\line(0,1){13}}
\put(13,0){\line(0,1){13}}
\put(0,0){\line(1,0){13}}
\put(0,4){\line(1,0){13}}
\put(0,8){\line(1,0){13}}
\put(0,12){\line(1,0){13}}
\put(0,.5){$F$}
\put(0,2.5){$\overline{l}$}
\put(0,4.5){$X$}
\put(0,6.5){$\overline{H}$}
\put(0,8.5){$Z$}
\put(0,10.5){$\tau$}
\put(12,12.5){$F$}
\put(9.5,12.5){$\overline{l}$}
\put(8,12.5){$X$}
\put(6,12.5){$\overline{H}$}
\put(3.5,12.5){$Z$}
\put(1.5,12.5){$\tau$}
\put(1.5,.5){$1$}
\put(1.5,2.5){$1$}
\put(1.5,4.5){$\overline{l}$}
\put(1.5,6.5){$2 \overline{l}$}
\put(1.5,8.5){$\overline{H}$}
\put(1.5,10.5){$\overline{H}$}
\put(3.5,.5){$0$}
\put(3.5,2.5){$1$}
\put(3.5,4.5){$F$}
\put(3.5,6.5){$2 \overline{l}$}
\put(3.5,8.5){$2X$}
\put(3.5,10.5){$\overline{H}$}
\put(5.5,4.5){$1$}
\put(5.5,6.5){$2$}
\put(5.5,8.5){$2\overline{l}$}
\put(5.5,10.5){$2\overline{l}$}
\put(7.5,4.5){$0$}
\put(7.5,6.5){$1$}
\put(7.5,8.5){$F$}
\put(7.5,10.5){$\overline{l}$}
\put(6.5,2){0}
\put(9.5,10.5){$1$}
\put(9.5,8.5){$1$}
\put(11.5,10.5){$1$}
\put(11.5,8.5){$0$}
\put(10.5,5.5){$0$}
\put(10.5,2){$0$}
\end{picture}
\end{center}
If $rank(\G)=4$, \ $\G$ is a cone with vertex a line $r$ over a smooth
quadric
surface
$Q\subset \Pin{3}.$ Following
\cite{Alibabaquad} let
$\tau$ be the tautological divisor on $\Gamma^*$ and
Let $C(W)$ denote the cone with vertex
$r$ over the cycle $W\subset \overline{Q}.$ Let
\begin{xalignat}{3}
Pic \ \overline{Q} : &= <\ell_1, \ell_2> & A_0(\overline{Q}) : &= <p> &Q
:
&= \sigma^{-1}(\overline{Q})  \notag \\
P_1 :&= \sigma^{-1}(C( \ell_1)) & p_1: &=\tau \cdot P_1 & \elonebar:&=
\tau \cdot
p_1 \notag \\
 P_2 :&= \sigma^{-1}(C( \ell_2)) & p_2: &=\tau \cdot P_2 & \eltwobar:&=
\tau \cdot
p2  \notag \\
F :&= \sigma^{-1}(C(p))& \ell:&= \tau \cdot F & &  \notag
\end{xalignat}
With the above notation it is:
\begin{xalignat}{3}
Pic (\Gamma^*) &= <\tau, P_1, P_2> & A_2(\Gamma^*) &= < Q,p_1, p_2, F> &
A_1(\Gamma^*) &= < \elonebar, \eltwobar, \ell> \notag
\end{xalignat}

It is also  $T=\tau-P_1-P_2$, $c_1(\Gs)=4\tau-T$,
$c_2(\Gs)=3Q+4p_1+4p_2.$
Because $S'$ is an effective cycle it is 
 $S'=\alpha Q+\beta p_1+\gamma p_2+\delta F$ with 
$\alpha\geq 0$,
$\alpha +\beta\geq 0$, $\alpha+ \gamma \ge 0$ and $deg(S)=\tau\cdot
\tau\cdot
S'=2\alpha+\beta+\gamma+\delta.$ 

With the above notation we have the
 following intersection table:
\vskip .5cm

\setlength{\unitlength}{ .7cm}
\begin{center}
\begin{picture}(11,11)
\multiput(0,0)(0,3){2}{\line(1,0){11}}
\multiput(0,7)(0,3){2}{\line(1,0){11}}
\multiput(1,0)(3,0){2}{\line(0,1){11}}
\multiput(8,0)(3,0){2}{\line(0,1){11}}
\multiput(1.2,0.2)(1,1){10}{$1$}
\multiput(2.2,0.2)(1,1){9}{$0$}
\multiput(1.2,1.2)(0,1){2}{$1$}
\multiput(2.2,2.2)(0,1){2}{$0$}
\multiput(3.2,3.2)(0,1){2}{$0$}
\multiput(7.2,7.2)(0,1){2}{$0$}
\multiput(6.2,6.2)(1,0){2}{$1$}
\multiput(4.2,4.2)(0,1){2}{$1$}
\multiput(.2,9.2)(1,1){2}{$\tau$}
\multiput(.2,8.2)(2,2){2}{$P_1$}
\multiput(.2,7.2)(3,3){2}{$P_2$}
\multiput(.2,6.2)(4,4){2}{$Q$}
\multiput(.2,5.2)(5,5){2}{$p_1$}
\multiput(.2,4.2)(6,6){2}{$p_2$}
\multiput(.2,3.2)(7,7){2}{$F$}
\multiput(.2,2.2)(8,8){2}{$\overline{\ell_1}$}
\multiput(.2,1.2)(9,9){2}{$\overline{\ell_2}$}
\multiput(.2,.2)(10,10){2}{$\ell$}
\multiput(1.2,8.2)(1,1){2}{$p_1$}
\multiput(1.2,7.2)(2,2){2}{$p_2$}
\multiput(1.2,6.2)(3,3){2}{$\star$}
\multiput(1.2,5.2)(4,4){2}{$\overline{\ell_1}$}
\multiput(1.2,4.2)(5,5){2}{$\overline{\ell_2}$}
\multiput(1.2,3.2)(1,1){3}{$\ell$}
\multiput(5.2,7.2)(1,1){3}{$\ell$}
\multiput(3.2,6.2)(1,1){2}{$\overline{\ell_1}$}
\multiput(2.2,5.2)(3,0){2}{$0$}
\multiput(2.2,6.2)(2,2){2}{$\overline{\ell_1}$}
\multiput(2.2,7.2)(1,1){2}{$F$}
\multiput(2.2,8.2)(3,0){2}{$0$}
\multiput(3.2,7.2)(3,0){2}{$0$}
\multiput(8.2,9.2)(1,0){2}{$1$}
\put(1.2,9.2){$Q$}
\put(4.2,6.2){$2$}
\put(5.2,6.2){$1$}
\put(8.2,8.2){$0$}
\put(.2,10.2){$\cdot$}
\end{picture}
\end{center}

where $\star = \overline{\ell_1} + \overline{\ell_2}$ and the empty
spaces
are intended to
be $0.$

%
%
\section{Surfaces in $\Pin{4}$}
Throughout this section $S$ will be a surface embedded in
 $\Pin{4}$ by a very ample
line bundle $L,$  non degenerate of degree $9,$ with sectional genus
$g.$
Surfaces of degree $9$ in $\Pin{4}$ have been completely classified by
Aure and Ranestad in \cite{au-ra}. The investigation of the
projective normality of such surfaces is essentially contained in
their work. For completeness we present the global picture in this
section.
\begin{theo}[ \cite{au-ra}]\label{class}
 Let $\sel$ be as above and $\chi=\chi(\calo_S)$ then $S$ is a regular
surface with $K^2=6\chi-5g+23$, where:
\begin{itemize}
\item[1)] $g=6$, $\chi=1$ and $S$ is
rational or the projection of an Enriques surface of degree $10$ in
$\Pin{5}$ with center of projection on the surface;

\item[2)] $g=7$ and $\chi =1$ and $S$ is rational, or $\chi =2$ and $S$
is a minimal elliptic surface;
\item[3)] $g=8$ and $\chi =2$ and $S$ is a
K3-surface with $5$ $(-1)$-lines, or $\chi =3$ and $S$ is a minimal
surface of general type;
\item[4)] $g=9$, $\chi=4$ and $S$ is linked $(3,4)$
to a cubic scroll;
\item[5)] $g=10$, $\chi=6$ and $S$ is a complete
intersection $(3,3)$;
 \item[6)] $g=12$, $\chi = 9$ and $S$ is linked to a
plane.
\end{itemize}

Moreover if $g\geq 7$ then  $S$ is contained in at least two
quartic surfaces.
\end{theo}

\begin{prop}$\sel$ as above is projectively normal if and only if:
\begin{itemize}
\item [a)] $g=8$, $\chi =2$ and $S$ is a K3-surface with $5$
$(-1)$-lines.
\item[b)] $g=9$, $\chi=4$ and $S$ is linked $(3,4)$ to a cubic scroll
\item[c)] $g=10$, $\chi=5$ and S is a complete intersection
(3.3);
\item[d)]  $g=12$, $\chi = 9$ and $S$ is linked to a plane.
\end{itemize}
\end{prop}
\begin{pf} Let us examine the surfaces in Theorem \ref{class};\\ Let
$C\in|L|$ be a
generic smooth element.  Since all the surfaces are regular we always
have
$h^0(\restrict{L}{C})=4$. Note that for $g\leq 9$
$d(2\restrict{L}{C})<2g-2$, then
$h^1(2\restrict{L}{C})=0$. If $g=6$ then $h^1(\restrict{L}{C})=0$ and
thus
$h^1(L)=0$. Because
$h^0(2\restrict{L}{C})=h^1(2\restrict{L}{C})+19-6=13$    the following
exact
sequence
\begin{equation}
\label{duelle}
0\lra L\lra 2L\lra 2\restrict{L}{C}\lra 0
\end{equation}
gives $h^0(2L)=h^0(L)+h^0(2\restrict{L}{C})=18.$   Thus the map
$H^0({\cal O}_{{\bf P}^4}(2))\lra H^0({\cal O_S}(2))$ cannot be
surjective,
being  $h^0({\cal O}
_{{\bf P}^4}(2))=15.$ This means that  $S$ is not $2$-normal and
therefore it is not
projectively normal. \\ If $g=7$ then $h^0(2\restrict{L}{C})=12$,
$h^1(\restrict{L}{C})=1$ and
$h^1(2L)=0$ by
\cite[2.10]{au-ra}. Therefore from \brref{duelle} and the regularity of
S it follows that
$16\leq h^0(\restrict{L}{C})\leq 17$, which implies that $S$ is not
projectively
normal, as above.

 If $g=8$ then $h^0(2\restrict{L}{C}) =11$, $h^1(L)\leq 1$ and
$h^1(2L)=0$ by
\cite[2.11]{au-ra}.  For degree reasons $S$ cannot be contained in any
quadric hypersurface being contained in at least a quartic. Therefore
$S$
is $2$-normal if and only if $H^0(L)=15$.
>From \brref{duelle} we get
$ h^0(2L)\leq 15,16$ respectively if $\chi=2,3$. If $\chi=3$ $S$ is not
projectively normal as above.

 If $\chi=2$ $S$ is $2$-normal and therefore projectively normal by
Lemma
\ref{2norm}.

If $g=9$ $\sel$ is projectively normal by linkage, see
\cite[2.13]{au-ra}.

If $g=10$ $S$ is a complete intersection and thus projectively normal.

 If $g=12$ $S$ is linked to a plane and therefore it is arithmetically
Cohen-Macaulay by linkage, which implies $S$ projectively normal.
\end{pf}


\section{Surfaces  embedded in $\Pin{N} ,$ $N\ge 5.$}
Let $ S$ be  a smooth  surface, let $L$ be a very ample line bundle on
$S$
  and let 
$${\cal S}_g =\{(S, L) \text{ as above} 
\ |  L^2 = 9,
h^0(L)\ge 6, g(S, L)=g \text{ and}  \sel  \text{ is not a scroll}\}.$$
If
$\sel\in{\cal S}_g$,
by Castelnuovo's Lemma 
$g\leq 7.$
 Let $\cal{S} =
\bigcup_{g=0}^7
\cal{S}_g.$
 In the following lemmata a few preliminary results are collected.

\begin{lemma}
Let $\sel\in{\cal S}$ and let $C\in|L|$ be a smooth generic element.
Then
\begin{itemize}
\item[a)] if $h^1(\restrict{L}{C})=0$ then
$g(L)\leq5$;
\item[b)] if $h^1(\restrict{L}{C})=1$ then $g(L)=7$ or $6$;
\item[c)] if  $h^1(\restrict{L}{C})\geq 2$ then $g(L)=7$;
\end{itemize}\label{h1LC}
\end{lemma}

\begin{pf}
  a) For $g\geq 1$ we have $h^0(\restrict{L}{C})\geq 3$. If
$h^0(\restrict{L}{C})=3, 4$ then
$h^0(L)\le 4, 5$ respectively,  which is a contradiction. Thus
$h^0(\restrict{L}{C})
\ge 5$ i.e.
$5\leq h^0(\restrict{L}{C})=10-g$, i.e.
$g\leq 5.$\\
b)  Since $h^0(K_C-\restrict{L}{C})=1$ then $d\leq 2g-2$. Moreover being
$d$ odd it
is $d\leq 2g-3$

c) If $h^1(\restrict{L}{C})\geq 2$ then $K_S|_C$ is a special divisor on
 $C$ and it contributes to $cl(C),$ thus
$cl(C)\geq 0$,
i.e \\ $0\leq d(K_C-\restrict{L}{C})-2h^0(K_C-\restrict{L}{C})+2\leq
2g-13$
, from which we get
$g\geq 7$.
\end{pf}

\begin{lemma}\label{nonP41}
 Let $\sel \in \cal S$ . Then $\sel$ is
projectively normal or
\begin{itemize}
\item[a)] $\sel \subset \Pin{5}$, $g=4$, $S$ is $\Pin{1}$-bundle over an
elliptic curve,
$e=-1$, $L=3C_0.$
\item[b)] $\sel \subset \Pin{5}$ is a conic bundle over an elliptic
curve,
$g=4$.
\item[c)] $\sel \subset \Pin{5}$ $g = 5,6.$
\end{itemize}
\end{lemma}
\begin{pf}
In our hypothesis it is $\Delta(S, L)=11-h^0(L)$. If codim$ (S) =1$ then
$S$ is
projectively normal. Because $S$ is not embedded in  $\Pin{4}$, we can
assume
codim$ (S) \ge 3$ i.e. $h^0(L) \geq 6$, i.e $\Delta(S)\leq 5$.
 First assume $\Delta < 5.$ It is $9=d\geq 2\Delta + 1$, then if
$g\geq\Delta$
$\sel$ is projectively normal by Theorem \ref{fujitatheo}.
Because $g=0$ implies $\Delta = 0$ it follows that $\sel$ is
projectively
normal if
$\Delta=0,1.$
Moreover if $g=1$ then $\sel$ is an elliptic scroll which is impossible.
Therefore $\sel$ is projectively
normal  if
$\Delta = 2.$

Let now $\Delta=3.$ By
\cite{Io1} it is $g=3$ and therefore $\sel$ is projectively normal

Let now $\Delta= 4.$ By \cite{Io4} Theorem 3, $\sel$ is projectively
normal
unless,
possibly, if it is a
scroll over a curve of genus $g=2,$  which is impossible.

Let now  $\Delta(S)=5$, i.e $h^0(L)=6$.
If  $g(L)=7$ then  $S$ is a
Castelnuovo Surface, see \cite{Har1}, and thus
 projectively normal.
 Let $g=2,3.$ Simple cohomological computations, using the
classification
given in \cite{Io1} show that there are no such surfaces in $\Pin{5}.$

Let $g=4$ and let $C \in |L|$ be a generic hyperplane section. By Lemma
\ref{h1LC}
it is $h^1(\restrict{L}{C})=0$ and therefore $h^0(\restrict{L}{C})=6$ by
Riemann Roch. This shows that
$q(S) \ne 0.$ By \cite{LiAq} and \cite{Io4} the only possible cases are
b)
and c) in the
statement .
\end{pf}

\subsection{SECTIONAL GENUS $g=4$ }
\label{genusfour}
In this subsection  the projective normality of pairs $\sel\in{
\cal S}_4$ is studied. By Lemma \ref{nonP41} and \cite{LiAq} we have to
investigate
the following cases :
\begin{itemize}
\item[Case 1.]  $\sel \subset \Pin{5}$ is a $\Pin{1}$-bundle over an
elliptic curve, $e=-1$, $L\simeq3C_0.$
\item[Case 2.] $\sel \subset \Pin{5}$ is the blow up $p: S \to X$ of a
$\Pin{1}$-bundle over an elliptic  curve at  $3$ points,  $e=-1$,
$L\equiv  2\frak{C_0}+2\frak{f}-\sum_1^3E_i$ where $\frak{C_0} =
p^*(C_0)$ and
$\frak{f}=p^*{f}.$
\item[Case 3.]$\sel \subset \Pin{5}$ is the blow up $p: S \to X$ of a
$\Pin{1}$-bundle over an elliptic  curve at  $3$ points,  $e=0$,
$L\equiv 2\frak{C_0}+3\frak{f}-\sum_1^3E_i$ where $\frak{C_0} =
p^*(C_0)$ and
$\frak{f}=p^*{f}.$
\end{itemize}
\begin{lemma}
Let $\sel \in \cal{S}_4$ be as in \brref{genusfour} Cases 1,2,3 above.
The following
are equivalent:
\begin{enumerate}
\item[i)]$\sel$ is  $2$- normal;
\item[ii)]$\cal{I}_S$ is $3$-regular;
\item[iii)] $\sel$ is projectively normal;
\item[iv)] $h^0(\cal{I}_S(2)) = 0,$ i.e. $S$ is not contained
in any quadric hypersurface.
\end{enumerate}
\label{g4noquad}
\end{lemma}
\begin{pf}
Assume $\sel$ is $2$-normal, i.e. $h^1(\cal{I}_S(2)) = 0.$ Since $p_g(S)
=
0$ and
$h^1(L)=0$ in all the cases under consideration, it is not hard to check
that $\cal{I}_S$
is
$3$-regular. By
\cite[pg. 99]{mu1} it follows that
$h^1(\cal{I}_S(k)) = 0$ for all $k \ge 2$ and thus $\sel$ is
projectively
normal.
Therefore i), ii) and iii) are equivalent.  Since $h^0(\oofp{5}{2}) =
21$ and
$h^0(\oof{S}{2}) = 21,$
$S$ is $2$- normal if and only if $h^0(\cal{I}_S(2)) = 0.$
\end{pf}
\begin{rem}
$\sel$ as in \brref{genusfour} Case 3 is a congruence of lines of
$\Pin{3}$ of bi-degree
$(3,6),$  cf.
\cite{gr} and thus not projectively normal by Lemma \ref{g4noquad}.
\end{rem}
\begin{rem}
$\sel$ as in \brref{genusfour} Case 1 was shown to be projectively
normal by
Homma \cite{Ho2}. Following an idea due to Sommese we offer below a
different
proof.
\end{rem}
\begin{prop}
\label{propdiandrew}
Let $\sel$ be as in \brref{genusfour} Case 1, then $\sel$ is
projectively
normal.
\end{prop}
\begin{pf}
Lemma \ref{g4noquad} shows that it is enough to show that $\sel$ is
$2$-normal.
Let $E$ be an elliptic curve with  fixed origin $O.$ It was shown in
\cite{be-so2} that
$S$ can be viewed as the quotient of $X=E\times E$ under the involution
$ \iota:
X\to X$ given by $ \iota(x,y) = (y,x).$ Let $q: X \to S$ be the quotient
map
and $p_i: X \to E,$ $i=1,2$ be the projections onto the factors. One can
see that $
q^*(C_0) = p_1^*(\oof{E}{P}) \tensor p_2^*(\oof{E}{P}),$ where $P$ is a
point on $E.$ Let
$L_i = p^*_i(\oof{E}{P}).$ It is
$$ H^0(X, q^*(L)) = H^0(X, q^*(3C_0)) = H^0(3L_1) \tensor H^0(3L_2),$$
and therefore $h^0(q^*(L)) = 9.$

Let $H^0(X, q^*L)^{\iota}=\{ \sigma \in H^0(X,q^*L) \text{ such that }
\sigma \iota =
\sigma \}$, i.e. the subspace of global holomorphic  sections of
$q^*(L)$
which are
invariant under
$\iota.$ Notice that there is a natural isomorphism  $H^0(S, L) \simeq
H^0(X,q^*(L))^{\iota}.$
Because $deg (3L_i) \ge 2g(E) + 1$ the map $H^0(3L_i) \tensor H^0(3L_i)
\to
H^0(2(3L_i))$ is surjective. It follows that  $$H^0(q^*(L)) \tensor
H^0(q^*(L))  \to
H^0(q^*(2L))$$ is surjective.

To conclude it is enough to show that
\begin{equation}
\label{2normiota}
H^0(q^*(L))^{\iota} \tensor H^0(q^*(L))^{\iota}  \to
H^0(q^*(2L))^{\iota}
\end{equation}
 is surjective.
Let $\alpha_1, \alpha _2, \alpha_3$ be a base for $H^0(E, \oof{E}{3P})$
and let
$a_i=p_1^*(\alpha_i)$ and $b_i=p_2^*(\alpha_i)$ be bases for $H^0(3L_1)$
and
$H^0(3L_2)$ respectively.
With the above notation $\{ a_i \tensor b_j\}$ form a base for
$H^0(q^*(L)).$

Notice that an element $a_i \tensor b_j$ with $i \ne j$ cannot be
$\iota$-invariant.
To see this assume $a_i \tensor b_j$  were $\iota$-invariant and let $z
\in
E$ be a point
such that
$\alpha_i(z) = 0$ and
$\alpha_j(z) \ne 0.$ Then
$(a_i \tensor b_j)(z,y) = (a_i \tensor b_j)(y,z)$ implies
$\alpha_i(y)\alpha_j(z) = 0$ for
every
$y,$ which is impossible.  Simple direct checks show that
$$
a_1 \tensor b_1, \ a_2 \tensor b_2, \ a_3 \tensor b_3, a_1 \tensor b_2 +
\
a_2 \tensor
b_1, a_1 \tensor b_3 + \ a_3 \tensor
b_1, a_2 \tensor b_3 + \ a_3 \tensor
b_2
$$
form a basis for $H^0(q^*L)^{\iota}$  and that $$t_1= a_1 \tensor
b_2,t_2=
a_1 \tensor
b_3,t_3= a_2
\tensor b_3$$
span a $3$-dimensional subspace $V$ such  that $ H^0(q^*(L)) =
H^0(q^*(L))^{\iota}
\oplus V.$
To conclude the proof of the surjectivity of (\ref{2normiota}) it is
enough
to show
that elements of the form $v \tensor s$ or $s \tensor v$ where $v \in V$
and $s \in
H^0(q^*L)^{\iota}$ , and elements of the form $v_1 \tensor v_2$ where
$v_i
\in V$,
cannot be $\iota$-invariant.
Let $v= p_1^*(\gamma_1) \tensor p_2^*(\gamma_2)$ and let $s=
p_1^*(\sigma_1)
\tensor p_2^*(\sigma_2).$  For every $(x,y) \in X$ it is
\begin{gather}
\label{onesidevs}
(v \tensor s)(x,y) = \gamma_1(x) \gamma_2(y)\sigma_1(x) \sigma_2(y) =
\gamma_1(x) \gamma_2(y) \sigma_1(y) \sigma_2(x)\\
\label{othersidevs}
(v \tensor s)(y,x) = \gamma_1(y) \gamma_2(x) \sigma_1(y) \sigma_2(x).
\end{gather}
Equating \brref{onesidevs} and \brref{othersidevs} shows that  $v$ is
invariant, a
contradiction. The same argument takes care of the case $s \tensor v.$

Let now $v_1 = \sum c_i t_i$  and $v_2 = \sum d_i t_i,$ where the
$t_i's$
are as above
and assume that $v_1 \tensor v_2$ is invariant.  Since $\oof{E}{3P}$ is
generated by
global sections there is a point $z \in E$ such that $\alpha_1(z) \ne
0,$ while
$\alpha_2(z) = \alpha_3(z) = 0.$ 

The fact that $(v_1 \tensor v_2)(z,y) =
(v_1 \tensor
v_2)(y, z)$ for every $y$ implies$$(c_1 + d_1)\alpha_1(z) \alpha_2(y) =
0$$
for every $y,$ which gives $c_1 = -d_1.$ Repeating the argument
permuting the
indices it follows that $v_2 = - v_1.$
It is then enough to show that $v_1 \tensor (-v_1)$ is not
$\iota$-invariant. If it
were it would follow that $(v_1(x, y) - v _1(y, x))(v_1(x, y) + v_1(y,
x))
= 0$ for all
$x,y
\in E.$ Because $v_1$ is not invariant and it is  not  zero everywhere,
this is a
contradiction.
\end{pf}

\begin{lemma}\label{tipidicurve}
Let $\sel$  be as in \brref{genusfour} Case 2 . Then
\begin{enumerate}
\item[i)] If $r$ is a line contained
in $S,$ then
$r=E_i$  or $r=\frak{f}-E_i$, $i=1,2,3.$
\item[ii)] If $C\subset S$ is a reduced irreducible  cubic with $C^2 \ge
0$
then $C$ is
a  curve whose numerical class  is $C \equiv \frak{C_0} + \frak{f} -
\sum_1^3 E_i$
\item[iii)] If $C\subset S$ is a reduced irreducible  quartic  with $C^2
\ge 0$  then $C$
is a  curve whose numerical class  is one of the following:
\begin{enumerate}
\item[a)] $C\equiv \frak{C_0}$
\item[b)] $C \equiv \frak{C_0} +\frak{f} - E_i -E_j $
\end{enumerate}
\end{enumerate}
\end{lemma}
\begin{pf}
 Let $r=a\frak{C_0}+b\frak{f}-\sum_1^3a_iE_i$ be a line in $S$. Then
$L \cdot r=4a+2b-\sum_1^3a_i=1$ and
$0=2g(r)=2+a(a-1)+2b(a-1)-\sum_1^3a_i(a_1-1)$. It follows that
\begin{gather}
\sum_1^3a_i^2=a^2+3a+2ab+1\geq\frac{1}{3}(4a+2b-1)^2 \notag \\
\text{i.e.} \ \ \ \ \ 4b(b-1)+13a(a-1)+10ab-4a-2\leq 0
\label{diseqretta}
\end{gather}

Since $r$ is an irreducible smooth curve either $r=E_i$ or $r$ is the
strict
transform of an irreducible curve on the $\Pin{1}$-bundle and  the
following cases
can occur:
\begin{itemize}
\item $a=0$ and $b=1$, that gives us the fibers through the points blown
 up, i.e. $r=\frak{f}-E_i$ for $i=1,2,3$.

\item $a=1$ and $b\geq 0$ , for which \brref{diseqretta} would imply
$b=0$,
$\sum_1^3a_i =3$ and $\sum a_i^2=5$, i.e. $r=\frak{C_0}-2E_i-E_j$ for
$i,j=1,2,3$. But this
would imply the existence  of an irreducible  curve in $|C_0|$ passing
through a
point
$P_i$ with  multiplicity $2,$ that would imply $C_0 \cdot f=2$, where
$f$
is the fiber
through $P_i$, which is a contradiction.

\item $a\geq 2$ and $b\geq -\frac{a}{2}$. Let $a=2+h$ with  $h\geq 0.$
>From \brref{diseqretta} it follows that
 $4b^2+16+8h^2+27h+5h(h+2b)+8(2b+h)\leq 0$ and therefore
$ 4b^2+8h^2+17h\leq 0$, which is impossible unless $b=0$, $a=2$ which
 contradicts \brref{diseqretta}.
\end{itemize}
Cases ii) and iii) follow from similar computations and Castelnuovo's
bound
on the
arithmetic genus of curves.
\end{pf}

%
%
Techniques found in \cite{Alibabaquad} and a detailed analysis of the
geometry of
hyperplane sections will be used to deal with Case 2.
\begin{prop}
\label{selincones} Let $\sel \in \cal{S}_4$ be as in Case  2. Then
$\sel$ is
 projectively normal.
\end{prop}
\begin{pf}
By contradiction assume $\sel$ is not projectively normal.  Lemma
\ref{g4noquad}
then implies that
$\sel$ must be  contained in a quadric hypersurface $\Gamma.$ From
\cite{ar-so2}
it follows that $\Gamma$ must be singular.
\begin{case}  $rk(\G)=5.$
\end{case} Let $P$ be the vertex of the quadric cone $\Gamma.$ Following
the
notation of subsection
\ref{qcones}  let $S'=\alpha \overline{H} +X.$ If $P$ is not contained
in
$S$, then $S=S',$
$deg(S)=2\alpha$, which is impossible. We can assume $P\in S$. Then
$S'=4\overline{H}+X$ because $deg(S')=\tau^2  S'$
 and $E^2=(\restrict{T}{S'})^2=-1.$ Moreover $c_1(\Gs)|_{S'} \cdot
c_1(S')=2(2L-E)(-K_S'-E)$, since
$\tau|_ {S'}=\sigma^*(L)$.
 Plugging the above obtained values  into \brref{DPF}  a contradiction
is
reached.

\begin{case} $rk(\G)=4$
\end{case}
 Let $r$ be the line vertex of the cone.
>From $deg(S')=9$ we have
\begin{gather}
9=2\alpha+\beta+\gamma+\delta
\end{gather}
If $r\subset S$ then by Lemma \ref{tipidicurve} $r=E_i$ or
$r=\frak{f}-E_i$. Notice that in this case $S'= S.$ Let $T|_{S'}=\lambda
r,$ since
$T|_{S'}\cdot\tau|_{S'}=
\delta$ and $(T|_{S'})^2=-\delta^2$ we have $\lambda=\delta$ and
$\beta+\gamma=\delta-
\delta^2.$ Moreover $(S')^2=2(\alpha +\beta)(\alpha +\gamma)+2\alpha
\gamma$,
$c_2(\Gs)|_{S'}=14\alpha+7(\beta +
\gamma)+3\delta$, $c_1(\Gs)|_{S'}=4\sigma^*(L)-\delta r$. But from
$9=2\alpha+\beta+\gamma+
\delta$ the only possible values are
$(\alpha,\delta,\beta+\gamma)=(3,1,2)$
which give a
contradiction in  \brref{DPF}.

If $S\cap r= \emptyset$ then $9=2\alpha$, since $T|_{S'}=0,$ which is a
contradiction.

If $S\cap r=\{P_1,...,P_k\}$, let $\mu_j$ be the multiplicity of
intersection at $P_j$ and
let $s=\sum\mu_j$. Then  $(T|_{S'})^2=-\sum \mu_j=-s.$ If any of the
$\mu_j's$ is
strictly greater then $1,$ $S'$ acquires a singularity of type $A_{\mu_j
-
1}$ at a
point of $\overline{E_j},$ where $\overline{E_j}$ are the
exceptional divisors of $\restrict{\sigma}{S'}.$ Notice that $(\tau
T)S'=0$
gives
$S'=\alpha Q+\beta p_1+\gamma p_2$, with $\beta+\gamma=s.$
 Moreover it is $\alpha\geq 2$ because $\alpha$ can be viewed as the
degree
of the
generically finite rational map $\psi: S\lra Q$ induced by the
projection
from the
vertex of
$\G$, where $S$ is not birational to $Q$. Thus the only possible values
are
$(\alpha, s)=(4,1),(3,3),(2,5)$.  Let us assume at first that $\mu_j =
1$
for all $j,$ so
that $S'$ is smooth. 

It is  $(S')^2=2\alpha(\alpha+h) +2\beta\gamma$,
$c_2(\Gs)|_{S'}=63$,
$c_1(\Gs)|_{S'}=4\sigma^*(L)-\sum\overline{E_j}.$ Using the admissible
values for
$\alpha$ and $s$ we get a contradiction in \brref{DPF}.
Therefore for at least one $j$ it is  $\mu_j \ge 2$  and $(\alpha,s) =
(4,1)$ does
not occur.
Let $\Pi\subset \Pin{3}$ be a general 2-plane tangent to $Q.$ Then $\Pi
\cap Q = \ell_1
\cup
\ell_2$ where $\ell_i$ is a line in $\Pi.$ Cutting  $S$  with
the hyperplane spanned by $\Pi$ and $r$ we get a degree nine divisor
$D\in |L|$
which must be reducible as $D=D_1 \cup D_2$ where $\psi(D_i) =\ell_i.$
Moving $\Pi$ along $\ell_i$ we can see that $D_j$ moves at least in a
pencil. Therefore $h^0(D_i)\ge 2$ for $i=1,2.$ Moreover the above
argument shows
that $D_i$ is spanned away from $S\cap r.$ In particular $D_i$ cannot
have
a fixed
component, therefore $D_i^2 \ge 0.$  Let
$d_i = L
\cdot D_i.$  Then
$(d_1,d_2) = (1,8), (2,7), (3,6),(4,5).$

Lemma \ref{tipidicurve} shows that $S$
contains only a finite number of lines, therefore the first case cannot
happen.

In the
second case, moving $\Pi$ along $\ell_2$, $S$ could be given a conic
bundle
structure
over $\Pin{1}$ which is not possible.

When $(d_1,d_2) = (3,6)$ notice that $D_1$ must be reduced and
irreducible
because $S$ contains only a finite number of lines. Therefore $D_1\equiv
\frak{C_0}
+
\frak{f} - \sum_i E_i$ as in Lemma
\ref{tipidicurve} ii). When $s=3$, $\psi $ is a generically $3:1$ map
while
when $s=5$
$\psi
$ is a generically $2:1$ map. Therefore there is always at least a point
$P\in S\cap
r$  such that
$P
\not \in D_1.$
pertanto
meno con
 Because $h^0(X, C_0 + f)=3$ and $C_0 + f$ is spanned,
$h^0(D_1) \le 2.$ Since $|D_1|$ must be at least a pencil, it is
$h^0(D_1)
= 2.$ This
shows that the complete linear system $|D_1|$ is obtained by moving
$\Pi$ along
$\ell_2.$  A member
of
$|D_1|$ passing through $P$
 can then be found, contradiction.

Let now $(d_1, d_2) = (4,5).$ Assume $D_1$ reduced and irreducible. Then
$D_1$
must be as in Lemma \ref{tipidicurve} iii). Because $h^0(\frak{C_0}) =1$
it is
$D_1\equiv \frak{C_0} + \frak{f} - E_i -E_j.$ We claim that $D_1$ is
then a
smooth
elliptic quartic embedded in
$\Pin{3}.$ To see this notice that every element of $|C_0 + f|$ on $X$
is smooth
with the only exception of one curve, reducible as the union $C_0 \cup
f.$
Moreover notice that the same argument used above shows that $h^0(D_1) =
2$
and $|D_1|$
 is obtained by moving $\Pi$ along $\ell_2.$  Because
$\psi(D_1) =\ell_1$, for degree reasons
$D_1$ must go through at least one point in $S\cap r.$ Because $\mu_j\ge
2$
for at
least one $j$ and $h^0(D_1)=2$ we can always assume that $D_1$ has a
$(k\ge
3)$-secant  line.

It is known (see \cite{Io2}) that the ideal of such quartics in
$\Pin{3}$
is generated
by quadrics and therefore they cannot have $(k \ge 3)$-secant lines.
Let $D_1$ now be reducible or non reduced. $D_1$ cannot be reducible
with
lines as
components since $S$ contains only a finite number of lines. A simple
numerical check shows that the only smooth conics on $S$  have numerical
class
$\frak{f}$. Therefore we can assume $D_1 \equiv 2\frak{f.}$ As it was
pointed out above $D_1$
must pass through at least a point $P \in S \cap r$ but this contradicts
$\frak{f}^2 =
0.$
\begin{case} $rk(\G)=3$
\end{case}
Assume $S \subset \Gamma$ where $\Gamma$ is a quadric cone with $rk 
\Gamma
= 3$ and vertex $ V \simeq \Pin{2}$ over a smooth conic
$\gamma\subset\Sigma
\simeq \Pin{2}.$   Let $\psi: S -->\gamma$ be the rational linear
projection from $V.$
Let $\ell_1, \ell_2$ and $\ell_3$  be distinct  lines in $\Sigma$ such
that
$\gamma \cap
\ell_1 =\{P_1, P_2\}$,\ \ $\gamma \cap \ell_2
=\{P_2, P_3\}$, and  $\gamma \cap \ell_3
=\{P_1, P_3\}$, where $P_i \neq P_j.$ Let $D_i$ be the hyperplane
sections
of $S$
given by the hyperplanes spanned by $V$ and $\ell_i.$ Notice that $D_i$
must be
reducible.

Assume that $D_i$ has no components contained in $V$ for at least one
$i$, say
$i=1.$  Then
$D_1\sim A+B$ where   $\psi(A)=P_1$ and $\psi(B)=P_2.$ Let $L\cdot A =
a$ and
$L \cdot B=b.$ It follows that $D_2 \sim A'+B$ and $D_3 \sim A'+A$.
Notice
that $A
\sim A'$ and so $L\cdot A'=a.$ This leads to the contradiction
$2a=L\cdot
D_3=9.$

It follows that for any hyperplane section $D$ obtained
with the hyperplane spanned by a line $\ell \subset \Sigma$ and $V$ it
must
be $D
\sim 2C + F$ where $F\subset V$ and no component of $C$ is contained in
$V.$
Because $|C|$ is at least a pencil and it  cannot clearly have fixed
components, it must
be $C^2\ge 0.$
Since $S$ contains only a finite number of lines and it is not a
rational
conic bundle
it is $L\cdot C = 3,4.$ If $L\cdot C =3$ $C$ must be irreducible and
therefore as in
Lemma \ref{tipidicurve}, i.e. $C\equiv \frak{C_0}+\frak{f} - \sum_i
E_i.$
It follows that
$F \equiv E_1+E_2+E_3$ which is impossible because these three lines are
disjoint.
If $L\cdot C=4$ and $C$ is reduced and irreducible then $C \equiv
\frak{C_0} +
\frak{f} - E_i-E_j$ as in Lemma \ref{tipidicurve} iii) (the case
$C\equiv
\frak{C_0}$ cannot happen since
$h^0(\frak{C_0})=1$). Then
$F \equiv E_i+E_j - E_k$ which is impossible. If $C$ is reducible or non
reduced then
$C \equiv 2\frak{f}$ and $F\equiv 2\frak{C_0}-2\frak{f} -\sum_iE_i $
which
is not
effective because
$C_0$ is ample on $X$ and $C_0\cdot(2C_0-2f) =0.$
\end{pf}

%
%
\subsection{SECTIONAL GENUS $g =5$}

In this section we will study the projective normality of pairs
$\sel\in{
\cal S}_5.$ From Lemma \ref{nonP41} it follows that either $\sel$ is
known to
be projectively normal or $S\subset
\Pin{5}$ and hence
$\Delta(S, L)=5.$
 Since $g=\Delta$ and $d\geq 2\Delta-1$ the ladder
is regular and therefore $\sel$ is projectively normal if
$(C,\restrict{L}{C})$ is
projectively
normal.
\begin{lemma}
\label{equivalent}
 Let $\sel\in{\cal S}_5$ then if $h^1(L)=0$  the following statements
are  equivalent:
\begin{itemize}
\item[1)] $\sel$ is projectively normal
\item[2)] $(C,\restrict{L}{C})$ is projectively normal
\item[3)] $S$ is contained in exactly one quadric hypersurface in
$\Pin{5}$.
\end{itemize}
\end{lemma}
\begin{pf} Since $h^1(L)=0,$ $\sel$ is projectively normal if and only
if
$(C, \restrict{L}{C})$
is projectively normal by Lemma \ref{besanaignorans}. Moreover $(C,
\restrict{L}{C})$ is
projectively normal if and only if it is 2-normal, by Lemma \ref{2norm}.
Consider the exact sequence:
$$0\lra{\cal I}_S(1)\lra{\cal I}_S(2)\lra{\cal I}_C(2)\lra 0.$$
It is $h^1({\cal I}_S(1))=0$ since $L$ is linearly normal and $h^0({\cal
I}_S(1))=0$ since
$S$ is non degenerate. It follows that $h^0({\cal I}_S(2))=h^0({\cal
I}_C(2)).$
Because $h^0({\cal O}_{\Pin{4}}(2))=15$ and $h^0(2\restrict{L}{C})=14$
we
have that
$(C, \restrict{L}{C})$ is
projectively normal if and only if $h^0({\cal I}_C(2))=1$ and therefore
$\sel$ is
projectively normal if and only if $S$ is contained in exactly one
quadric.
\end{pf}
\begin{lemma} \label{trighyper}
Let $\sel \in \cal{S}_5.$ Then either $\sel$ is projectively normal or
it has a
hyperelliptic or trigonal section $C \in |L|.$
\end{lemma}
\begin{pf}
Since
$g(C)=5$ then
$cl(C)\le 2.$ If $cl(C)=2$ then by Theorem \ref{glcliff} $(C,
\restrict{L}{C})$ is
projectively normal and thus $\sel$ is projectively normal by regularity
of
the ladder.
\end{pf}
\begin{lemma}
\label{cbhyper} Let  $\sel\in{\cal S_5}$ with  a hyperelliptic section
$C\in |L|.$ Then $S$ is
a  rational conic bundle, not projectively normal.
\end{lemma}
\begin{pf} Surfaces with hyperelliptic sections are classified in
\cite{so-v}.
By degree considerations the only possible case is a rational conic
bundle.
Since
$h^1(L)=0$
 Proposition \ref{hyper} and Lemma \ref{equivalent} imply that $\sel$ is
not
projectively normal.
\end{pf}
\begin{lemma}
\label{cbtrig}
There are no  conic bundles $\sel$  with a trigonal
section $C\in |L|$ in ${\cal S}_5$ .
\end{lemma}
\begin{pf}  If $S$ is a conic bundle with trigonal section then
\cite{fa} Lemma 1.1 gives 
$g=2q+2$  which is impossible because \cite{LiAq} gives  $q\leq 1.$
\end{pf}

\begin{theo}
\label{g5theo}
 Let $\sel\in{\cal S}_5.$ Then $\sel$  fails to be projectively
normal if and
only if it is
\begin{itemize}
\item[1)] A  rational conic bundle;
 \item[2)] $(S, L)=(Bl_{12}\bold{ F_1}, 3\frak{C_0}-5\frak{f}-12p)$ with
trigonal section $C \in |L|$
and \\$\restrict{L}{C}=K_C-g^1_3+D$.
\end{itemize}
\end{theo}
\begin{pf}
>From \cite{LiAq}, Lemma \ref{trighyper}, Lemma
\ref{cbhyper}
and Lemma \ref{cbtrig}, the following cases are left to investigate:
\begin{center}
\begin{tabular}{|l|l|l|l|r|} \hline
Case &$S$ & $L$& existence\\ \hline\hline
1 & $Bl_{10}\Pin{2}$ & $7p^*(\cal O_{\Pin{2}}(1))-10\sum_1^{10} 2E_i$ &
Yes \\ \hline
2 & $Bl_{12}\Pin{2}$ & $6p^*(\cal O_{\Pin{2}}(1))-\sum_1^5
2E_i-\sum_6^{12} E_j$ &Yes\\ \hline
3 & $Bl_{10}\bold {F_e}$, $e=0,1,2$ &
$4\frak{C_0}+(2e+5)\frak{f}-\sum_1^7
2E_i-\sum_8^{
10} E_j$ & Yes\\ \hline
4 & $Bl_{12} \bold{ F_1}$ & $3\frak{C_0}+5\frak{f}-\sum_1^{12}E_i$ & ?
\\ \hline
\end{tabular}
\end{center}
where the hyperplane section is trigonal.
In case 2 $\sel$ admits a first reduction $(S',L')$ with $d'=16$. By
\cite{fa} $\sel$ cannot have trigonal section and therefore it is
projectively normal.

In case 3 it is  $K_S^2=-2$ so that $K_S(K_S+L)=-3$. Then by
\cite{bri-la}
Theorem 2.1
$\sel$ cannot have trigonal section thus it is projectively normal.

Case 1 is a congruence of lines of $\Pin{3}$ of bi-degree $(3,6)$
studied in detail
in \cite{ar-so2}.
In particular if ${\cal I}_S^*$ is the ideal of $S$ in the grassmanian
$G(1,3)$ of
lines of $\Pin{3},$ it is $h^0({\cal I}_S^*(2))=0$. From:
$$ 0\lra {\cal I}_G\lra{\cal I}_S\lra{\cal I}_S^*\lra 0$$
recalling that $G\in|{\cal O}_{\Pin{5}}(2)|$ we get $h^0({\cal
I}_S(2))=1$
and therefore
$\sel$ is projectively normal by Lemma \ref{equivalent}.

By \cite{GL} Corollary 1.6  a trigonal curve of genus 5 and degree 9 in
$\Pin{4}$ fails to
be projectively  normal if and only if it is embedded via a line bundle
$\restrict{L}{C}=K_C-g^1_3+D$ where
$D$ is an effective divisor of degree 4. Notice that this means that $
C$
is embedded in
$\Pin{4}$ with a foursecant line.
\end{pf}

\begin{rem} (ADDED IN PROOF)
\label{bali}
After this work was completed the first author and A. Alzati proved in
in \cite{bali}
that there exist  no surfaces  as in Theorem  \ref{g5theo} Case 2).
Therefore
this  case does not appear in the table of Theorem 1.1.
\end{rem}
%
%
\subsection{SECTIONAL GENUS $g= 6$}
\label{genere6}

In this subsection we will study the projective normality of pairs
$\sel\in{
\cal S}_6$. Notice that by Lemma \ref{nonP41}  $L$ embeds $S$ in
$\Pin{5}$
and the
ladder is regular.
\begin{lemma}
\label{2normg6}Let $C$ be a curve of genus $6$ embedded in $\Pin{4}$ by
the
complete linear system associated with  a very ample line bundle
$\restrict{L}{C}$ of
degree
$9$. Then
$C$ is
$2$-normal.
\end{lemma}
\begin{pf} Consider the exact sequence: $$0\lra {\cal I}_C\otimes{\cal
O}_{
\Pin{4}}(2)\to{\cal O}_{\Pin{4}}(2)\to{\cal O}_C(2)\lra 0$$
 Since $h^0({\cal O}_{\Pin{4}}(2))=15$ and $h^0({\cal
O}_C(2))=h^0(2\restrict{L}{C})=18+1-6=13$ it is $h^0(\iof{C}{2}) \ge 2$
with  the
map
$H^0({\cal O}_{\Pin{4}}(2))\lra H^0({\cal O}_C(2))$  surjective if and
only
if  $h^0(
{\cal I}_C(2))=2.$  Assume $h^0(\iof{C}{2})\ge 3,$ and let $Q_i$ for $i
=1,2,3$ be
three linearly independent quadric hypersurfaces containing $C.$ Because
deg$C
=9$  and $C$ is non degenerate it must be $dim(\cap_iQ_i) = 2.$ Let
$\frak{S} =
\cap_iQ_i$ then $h^0(\frak{S}) \ge 3.$ But $\frak{S}$ is a complete
intersection
$(2,2)$ in $\Pin{4}$ and it is easy to see that $h^0(\frak{S}) = 2$,
contradiction.
\end{pf}
\begin{cor}\label{gen6} Let $\sel \in \cal{S}_6.$ Then $(S,L)$ is
projectively normal.
\end{cor}
\begin{pf} Let $C\in|L|$ be a generic section. From Lemma \ref{2normg6}
and
Lemma
\ref{2norm} it  follows that $(C,\restrict{L}{C})$ is projectively
normal.
Since the
ladder is regular
this implies that $(S, L)$ is projectively normal.
\end{pf}
\section{Results on Scrolls}
A $n$-dimensional polarized variety $\xel$ is said to be a scroll over a
smooth curve
$C$ of genus $g$ if there is a vector bundle  \map{\pi}{E}{C} of rank $r
=
rk\,E =
n+1$  such that $\xel \simeq ({\Bbb P}(E), \oof{{\Bbb P}(E)}{1}).$

Recall that given a vector bundle $E$ over a curve $C$,
$\mu (E)$ and $\mu^- (E)$ of
$E$ are defined as ( see \cite{bu} for details) $$\mu(E) = \frac{deg
E}{rk E} =
\frac{d}{r}.$$
$$\mu^-(E)=min\{\mu(Q)|E\rightarrow Q\rightarrow 0\}$$

 $E$   will be called {\em very
ample} to signify that the tautological line bundle $\taut{E}$ is a very
ample line bundle on
${\Bbb P}(E).$
From \cite[Th 5.1.A]{bu} and general properties of projectivized bundles
it follows that:
\begin{prop}[\cite{bu}]\label{buscroll}  Let $\xel$ be a scroll. If
$\mu^-(E)>2g$ then
$\xel$ is projectively normal.
\end{prop}
The following Lemma is essentially due to Ionescu \cite{fa-li9} :
\begin{lemma}[Ionescu]\label{scrolls}
Let $(X,L)$ be an n-dimensional scroll over a hyperelliptic curve $C$ of
genus $g$
with $L$ very ample. Then
 $\Delta=ng.$
\end{lemma}
\begin{pf}
a) Let $X=\Proj{E}$, $L=\taut{E}$ and $\pi:X\lra
C$.\\ By the Riemann-Roch theorem and the fact that $\pi_*({\cal
O}_X(1))=E$ it follows that:
$$h^0(L)=h^0(C,E)=h^1(C,E)+d-n(g-1)$$
Thus it is enough to show that $h^1(C,E)=0$.\\
By Serre duality $h^1(C,E)=h^0(C,K_C\bigotimes E^{*})$. Assume
$h^1(C,E)\neq 0.$ A non trivial section  $\sigma\in
h^0(C, K_C\bigotimes E^*)$ gives the following
surjection:
\begin{equation}
\label{sigmasur}
 (K_C\bigotimes E^*)^*\lra{\cal O}_C(-D)\lra 0
\end{equation} where
$D$ is the divisor on $C$ associated to $\sigma$. Tensoring
\brref{sigmasur} with $K_C$ we obtain
$$E\lra K_C-D\lra 0.$$ Because $\taut{E}$ is very ample,
$K_C-D$ is very ample on C. Moreover $K_C-D$ is a special line bundle on
$C$ because $h^1(K_C-D)=h^0(D)>0$. This is impossible because  $C$ is
hyperelliptic.
\end{pf}
\begin{lemma}
\label{scrollg2}
Let $(S,L)$ be a  two-dimensional scroll over a curve of genus 2 and
degree 9 in
${\Bbb P}^{6}$ with
$L$ very ample. Then
$X={\Bbb P}(E)$ with $E$ stable and $\xel$ is projectively normal.
\end{lemma}
\begin{pf} Let $S={\Bbb P}(E)$ where $E$ is a rank $2$ vector bundle of
degree $9$ over a smooth curve
 of genus $2$. If $E$ is stable then $\mu^-(E)=\mu(E)=\frac{9}{2}>4$ and
so by Proposition \ref{buscroll} $(S,L)=({\Bbb P}(E),{\cal O}_{{\Bbb
P}(E)})$ is projectively normal.

Assume now $E$ non stable. Then there exists a line bundle $Q$ with
deg($Q)\leq 4$ such that $E\to Q\to 0$. This contradicts the very
ampleness
of $Q$
as a quotient of  a very ample $E.$
\end{pf}
\begin{prop}
\label{scrollprop}
Let $\sel$ be a scroll of degree $d=9$ over a smooth curve $C$ of genus
$g.$  Then
$\sel$ is projectively normal unless possibly if $C$ is trigonal, $3\le
g
\le 5$ and
$S\subset \Pin{5}.$
\end{prop}
\begin{pf}
Following the proof  of Lemma \ref{nonP41} if $\Delta \ge 2$ and  $g=1$
then
$\sel$ is an elliptic scroll (see \cite{fu})
and therefore projectively normal by
\cite{Alibaba} or \cite{Ho1},\cite{Ho2}.
If $\Delta =4$ and $g=2$ then $\sel $ is projectively normal by Lemma
\ref{scrollg2}.
Let $\Delta=5.$ If $g=6$ then $\sel$ is projectively normal by
\ref{gen6}.
If $g=5$
by Theorem \ref{fujitatheo}, \ref{fujitatheo2}, \ref{glcliff}
$\sel$ is projectively normal unless $cl(C)\le 1.$ If $g=3,4$ then it is
always $cl(C)\le 1.$
By Lemma \ref{scrolls} $C$ must be trigonal.
\end{pf}
\section{An Adjunction Theoretic Problem}
\label{K+L}
The question of finding examples for the problem posed by Andreatta, Ein
and
Lazarsfeld (see introduction) is addressed below.
\begin{cor}
Let $(S, L)$  be a surface polarized with a very ample line bundle of
degree $d=9$
such that the embedding given by $|L|$ is not projectively normal. Then
there
does not  exist a very ample line bundle $\cal{L}$ such that $L= K_S +
\cal{L}$ unless
$(S, L)$ is the blow up of an elliptic $\Pin{1}$-bundle  as in the first
case of Theorem
\ref{thetheorem}.
\end{cor}
\begin{pf}
Let $\sel$ be as in the Table of Theorem \ref{thetheorem}, not as in the
first case. Assume $L=K +
\cal{L}$ with $\cal{L}$ very ample. Computing  $\cal{L}^2$ and
$g(\cal{L})$ and
using \cite{LiAq} lead to a contradiction
in every case.
Similarly a contradiction is reached if $\sel$ is a scroll over a curve
of
genus $3,4,5.$
\end{pf}
\begin{rem}
The existence of an example  of a surface as in case 1 of Theorem
\ref{thetheorem}
where
$L=K+\cal{L}$ with $\cal{L}$ very ample is a very delicate question.
Let $E$ be an indecomposable rank $2$ vector bundle over an elliptic
curve
$\cal{C}$ with
$c_1(E) =0$  and let
$X =\Proj{E}.$ Let $C_0$ be the fundamental section, let $M$ be any line
bundle
whose numerical class is $2C_0 + f$ and  let $p: S=Bl_3X
\to X$ be the blow up
of $X$ at three points $P_i$ $i=1,..,3.$ Using the same notation for the
blow up
introduced in  subsection
\ref{notation} consider a line bundle 
$L\equiv 2\frak{C_0} + 3\frak{f} -
\sum_iE_i.$
Notice that $L\equiv  K_S + \cal{L}$ where $\cal {L}\equiv 4\frak{C_0} +
3\frak{f} -
\sum_i2E_i.$
Moreover $\cal{L} \equiv K_S + H,$  $H^2 = 9,$  $H \equiv 3T$ where $T=
p^*(M) -
\sum_1^3 E_i.$  Recent results of Yokoyama and Fujita \cite{fuyo},
\cite{Yoko})  show  that the $P_i's $ can be chosen generally enough to
have 
$T$  ample but not effective.
Reider's theorem then  shows that
$\cal{L}$ is very ample if it is possible to choose the $P_i's$ such
that  for
\underline{every}  line bundle $M$ whose numerical class is $2C_0 + f$
it is
$|p^*(M)-\sum_i E_i|= \emptyset.$
\end{rem}


\end{document}